\documentclass[letterpaper,10 pt, conference]{ieeeconf}   

\IEEEoverridecommandlockouts                              
\usepackage{graphics} 
\usepackage{epsfig} 
\usepackage{amsmath} 
\usepackage{amssymb}  
\usepackage[latin1]{inputenc}

\newcommand{\beq}{\begin{equation}}
\newcommand{\eeq}{\end{equation}}
\newcommand{\bea}{\begin{eqnarray}}
\newcommand{\eea}{\end{eqnarray}}
\newcommand{\bef}{\begin{figure}}
\newcommand{\eef}{\end{figure}}
\newcommand{\bsc}{\begin{scriptstyle}}
\newcommand{\esc}{\end{scriptstyle}}
\newcommand{\bd}{\begin{displaymath}}
\newcommand{\ed}{\end{displaymath}}
\newcommand{\nn}{\nonumber}
\newcommand{\nnl}{\nonumber \\}
\newtheorem{lemma}{Lemma}[section]
\newtheorem{proposition}{Proposition}[section]

\def\bmat{\left[ \begin{array}}
\def\emat{\end{array} \right]}
\newcommand{\defn}{\stackrel{\triangle}{=}}
\renewcommand{\thesubsection}{\Alph{subsection}}
\def\qed{\hfill \vrule height 7pt width 7pt depth 0pt \smallskip}

\title{\LARGE \bf On the convergence of a Risk Sensitive like Filter}

\author{Mattia~Zorzi, Bernard C.~Levy\thanks{This work has been partially supported by the FIRB project ``Learning
meets time'' (RBFR12M3AC) funded by MIUR.}  \thanks{M. Zorzi is is with the
Dipartimento di Ingegneria dell'Informazione, Universit\`a di
Padova, via Gradenigo 6/B, 35131 Padova, Italy,
({\tt\small zorzimat@dei.unipd.it})} \thanks{B. Levy is with the Department of Electrical and Computer
Engineering, 1 Shields Avenue, University of California, Davis,
CA 95616 ({\tt\small bclevy@ucdavis.edu})}}

\begin{document}

\maketitle
\thispagestyle{empty}
\pagestyle{empty}

\begin{abstract}
In this paper, we analyze the convergence of a risk sensitive like filter where the risk sensitivity parameter is time varying. Such filter has a Kalman like structure and its gain matrix is updated according to a distorted version of the Riccati iteration. We show that the iteration converges to a fixed point by using the contraction analysis.
\end{abstract}

\section{Introduction}
Physical systems are often modeled by (nominal) linear models. One reason is that the corresponding  filtering problem is tractable, for instance if we consider the Gauss-Markov state space model then we obtain the Kalman filter. On the other hand, linear models are rather simple and thus introduce modeling errors. This implies that
the optimal filter may not perform well in the realty.

One possible strategy to deal with such problem is to use robust filtering. The pioneering works are due by Kassam, Poor and their collaborators, \cite{KASSAM},\cite{POOR}. This paradigm can be sketched as follows. One player (say, nature) selects the least favorable model in an allowable neighborhood about the nominal model, while the other player designs the optimal filter based on that least favorable model. 
Therefore, the optimal filter is obtained by solving a minimax problem. However, the implementation of such a filter can be very difficult because it depends on the characterization of the allowable neighborhood. To overcome this difficulty, a new class of robust filters based on the minimization of risk sensitive functions, which penalize large estimation errors, was introduced in  \cite{SDJ},\cite{Whi},\cite{BS}. The sensitivity to large errors is tuned by a risk sensitivity parameter. This approach considers Gauss-Markov state space models and the resulting robust filter is a Kalman like filter where the gain matrix is updated according to a distorted version of the Riccati iteration (say, risk sensitive Riccati iteration). Unfortunately, this method only considers  the nominal model. In \cite{LN}, a new minimax robust state space filtering problem was examineted. In this approach, at each time step all possible increments of the state space model are described by a ball about
the nominal increment. Its radius is fixed {\em a priori} and represents the tolerance budget available
at each time step. Therefore, the nature selects the least favorable model increment in the allowable ball, and the other player designs the optimal filter based on that least favorable model. It turns out the resulting robust filter is a risk sensitive like filter where the risk sensitivity parameter is now time varying. Accordingly, the gain matrix updating is governed by a risk sensitive like Riccati iteration.

An important issue for Kalman like filters is their convergence. In \cite{Bou}, under the assumption that the Gauss-Markov state space model is reachable and observable, it has been
established that the Riccati mapping is a contraction for the Riemann metric associated to the cone of positive definite matrices, and thus the Riccati iteration asymptotically converges. The same result can be proved by using the Thompson part metric \cite{LeL,GAUBERT_2012}. In \cite{LZ}, a similar contraction analysis has been considered to prove the convergence of the risk sensitive Riccati iteration. Here, the problem has been formulated in Krein space, see \cite{HSK2}, \cite{HSK3}. Then, it has been  shown that the $N$-block risk sensitive Riccati mapping is strictly contractive for the Riemann metric by choosing the risk sensitivity parameter sufficiently small. Regarding the risk sensitive like Riccati iteration, it seems to converge \cite{LN},  but no convergence result has been proved yet.

In this paper, we consider a similar contraction analysis to prove that  the risk sensitive like filter in \cite{LN} asymptotically converges for tolerance values sufficiently small. More precisely, we formulate the filtering problem in Krein space and then we show that the N-block risk sensitive like Riccati mapping is strictly contractive
provided that the time varying risk sensitivity parameter is smaller than a constant parameter. Moreover, it is always possible to find a lower bound of this iteration after a finite number of steps. As we will see, both the constant parameter and the lower bound allow to characterize a range of values of the tolerance for which the iteration converges.

\vspace{-0.05cm}
The paper is organized as follows.
 In Section \ref{sec:robust_filtering} we review the risk sensitive like filter presented in \cite{LN}.
 In Section \ref{sec:riemann} we review the Thompson part metric and contractive mappings needed for the following contraction analysis. In Section \ref{sec:Nblock} we construct the $N$-block risk sensitive like Riccati mapping. In Section \ref{sec:contraction_property} we characterize a range of values of $c$ for which the mapping is a strict contraction. Finally, in Section \ref{sec:example} we provide an example.
Throughout the paper, $\cal P$ denotes the cone of positive definite symmetric matrices, and $\overline{\cal P}$ its closure. 
Given $P\in\cal P $, $\lambda_1(P)\geq \lambda_2(P)\geq \ldots \lambda_n(P)>0$ are its eigenvalues
sorted in decreasing order.

\section{Risk sensitive like Filtering}
\label{sec:robust_filtering}

Consider a discrete-time stochastic process $y_t$ described by a nominal Gauss-Markov state space model
of the form \bea
x_{t+1} &=&A x_t + B u_t \label{3.1} \\
y_t &=& C x_t + D v_t, \;\; t\geq 0  \label{3.2}
\eea
where the state $x_t \in \mathbb{R}^n$, the process noise $u_t\in\mathbb{R}^{m}$, and the observation noise $v_t\in\mathbb{R}^p$.
The noises $u_t$ and $v_t$ are assumed to be zero-mean WGN processes
with normalized covariance matrices and independent, that is
\[
E \Big[ \left[
          \begin{array}{c}
            u_t \\
            v_t \\
          \end{array}
        \right]\left[
                 \begin{array}{cc}
                   u_s & v_s \\
                 \end{array}
               \right]
 \Big] = \left[
           \begin{array}{cc}
             I_m & 0 \\
             0 & I_p \\
           \end{array}
         \right] \delta_{t-s} \: ,
\]
where $\delta_t$ denotes the Kronecker delta function. The initial
state vector $x_0$ is assumed independent of the noises
$u_t$ and $v_t$ with nominal probability density
\[ p_0(x_0)\sim {\cal N} (0,P_0).\]
The pairs $(A,B)$ and $(A,C)$ are assumed to be reachable and observable, respectively. Moreover, we assume that the noises $u_t$ and $v_t$
affect all the components of the dynamics (\ref{3.1}) and observations (\ref{3.2}), that is $BB^T$ and $DD^T$ are positive definite. As observed in \cite{LN}, this is a natural
property to demand when the relative entropy is used to measure the proximity of statistical models, see below.

The robust filter proposed in \cite{LN} is designed according to the minimax point of view. More precisely, at time $t$, 
the nature selects the least favorable increment of the state space model in a ball about the nominal increment given by  (\ref{3.1})-(\ref{3.2}). Such a ball is characterized by requiring that the {\em Kullback-Leibler} divergence,
\cite{COVER_THOMAS}, between the two model increments is smaller than or equal to the tolerance $c\geq 0$. Note that, $c$ is fixed by the user. More precisely, the larger $c$ is, the worse increments the nature can select.

It turns out that the robust estimator $\hat x_{t+1}$ of $x_{t+1}$, given the observations $y_t, y_{t-1},\ldots ,y_0$,  
 obeys the Kalman like recursion
\beq
\hat{x}_{t+1} = A \hat{x}_t + K_t \nu_t
\: , \label{3.3}
\eeq
where
\beq
\nu_t \defn y_t - C\hat{x}_t  \label{3.4}
\eeq is the innovations process.
In (\ref{3.3}), the gain matrix
\beq
K_t = A(P_t^{-1}-\theta_{t-1} I_n)^{-1}C^T(R_t^{\nu})^{-1} \: , \label{3.5}
\eeq
where
\beq
R_t^{\nu} = E[\nu_t \nu_t^T] = C(P_t^{-1}-\theta_{t-1} I_n)^{-1}C^T +DD^T
\label{3.6}
\eeq
represents the variance of the innovations process, $\theta_{t-1}$ with \beq \label{interval_for_theta} 0<\theta_{t-1}< (\lambda_1(P_t))^{-1}\eeq is the unique solution to the equation
\beq \label{condizione_su_gamma}\gamma(\theta_{t-1},P_t)=c\eeq
where
\beq  \gamma(\theta,P)\triangleq \frac{1}{2}\left [\log\det(I-\theta P)+\mathrm{tr}[(I-\theta P)^{-1}]-n\right], \eeq
and if $\tilde{x}_t = x_t - \hat{x}_t$ denotes the
state prediction error, its variance matrix $P_t =
E[\tilde{x}_t \tilde{x}_t^T]$ obeys the distorted Riccati
iteration
\beq
P_{t+1}=  r^R_c(P_t) 
\label{3.7}
\eeq with initial condition $P_0$. The mapping $r^R_c(\cdot)$ is defined by
\beq r^R_c(P_t) \defn A[ P_t^{-1} + C^T(DD^T)^{-1}C -\theta_{t-1} I]^{-1} A^T  + BB^T.\nn \eeq

Note that, the robust filter (\ref{3.3})-(\ref{condizione_su_gamma}) is a risk sensitive like filter, \cite{Whi}, \cite[Chapter 10]{speyer_2008}.
In the classic formulation, however, the risk sensitive Riccati mapping is defined as
\beq
r^{RS}_\theta(P) \defn A[ P^{-1} + C^T(DD^T)^{-1}C -\theta I]^{-1} A^T + BB^T\nn
\eeq where the risk sensitivity parameter $\theta\geq 0$ is constant and does not depend on $P$. Moreover, for $\theta=0$ (risk neutral case) we obtain the Riccati mapping
\beq
r(P) \defn A[ P^{-1} + C^T(DD^T)^{-1}C]^{-1} A^T + BB^T \nn.
\eeq Finally, it is worth noting that (\ref{interval_for_theta}) implies that $r^R_c(P)\in{\cal P}$ for each $P\in{\cal P}$, that is $r_c^R(\cdot)$ is a mapping of ${\cal P}$. Such a property does not hold for the classic risk sensitive mapping because it may occur that $r^{RS}_{\theta}(P)\notin{\cal P}$ even when $P\in{\cal P}$,
\cite{LZ}.

\section{Thompson part metric and contraction mappings}\label{sec:riemann}
If $P$ is an element of ${\cal P}$
with eigendecomposition
\beq
P = U \Lambda U^T \label{2.1}
\eeq
where $U$ is an orthogonal matrix formed by normalized
eigenvectors of $P$ and $\Lambda = \mbox{diag} \, \{
\lambda_1 , \ldots , \lambda_n \}$ is the diagonal
eigenvalue matrix of $P$, the square-root
of $P$ is defined as
\[
P^{1/2} = U \Lambda^{1/2} U^T
\]
where $\Lambda^{1/2}$ is diagonal, with entries $\lambda_i^{1/2}$
for $1 \leq i \leq n$. Similarly, the logarithm of $P$ is
the matrix specified by
\[
\log(P) = U \log(\Lambda ) U^T \: ,
\]
where $\log(\Lambda)$ is diagonal with entries $\log(\lambda_i)$
for $1 \leq i \leq n$. Let $P$ and $Q$ be two positive definite
matrices of ${\cal P}$. Then $P^{-1}Q$ is similar to $P^{-1/2}QP^{-1/2}$,
so they have the same eigenvalues, and $P^{-1/2}QP^{-1/2}$ is
positive definite.
The Thompson part metric between $P$ and $Q$ is defined as
 \begin{eqnarray}
d_T(P,Q) &=& ||\log(P^{-1/2}QP^{-1/2})||_2\nn\\ &=& \max\{\log(\lambda_1(P^{-1}Q)),\log(\lambda_1(Q^{-1}P))\} \: ,\nn
 \end{eqnarray}
where $||\cdot ||_2$ denotes the spectral norm.

Let $f(\cdot)$ be a non expansive mapping of $\mathcal{P} $.
Its contraction coefficient (or Lipschitz constant) is defined as
\beq
\xi(f) = \sup_{P,Q \in {\cal P}, P \neq Q} \frac{d_T(f(P),f(Q))}{d_T(P,Q)}
\: . \label{2.5}
\eeq
From (\ref{2.5}) we get
\[
d_T(f(P),f(Q)) \leq \xi(f) d_T(P,Q).
\]  Moreover, if $\xi(f)< 1$, then $f$ is a strictly contractive mapping.

If
$f$ is a strict contraction of ${\cal P}$ for the metric $d_T$,
by the Banach fixed point theorem \cite[p. 244]{AE}, there
exists a unique fixed point $P$ of $f$ in $\bar{\cal P}$
satisfying $P=f(P)$. Moreover,
if the $N$-fold composition $f^N$ of a non-expansive mapping
$f$ (or simply $N$-block mapping of $f$) is strictly contractive, then $f$ has a unique fixed
point. Furthermore this fixed point can be
evaluated by performing the iteration $P_{k+1} = f(P_k)$
starting from any initial point $P_0$ of ${\cal P}$. We 
will consider in particular the Riccati like mapping
defined by
\beq 
f(P) = M[P^{-1} + \Omega ]^{-1}M^T + W \label{2.6}
\eeq
where $P$, $\Omega$ and $W$ are symmetric real positive definite
matrices and $M$ is a square real, but not necessarily invertible,
matrix. For this mapping the following result was established
in \cite[Th. 5.3]{LeL}.
\vskip 2ex

\begin{lemma}
\label{lem1}
The mapping in (\ref{2.6}) is strictly contractive with
\beq
\xi(f) \leq \left(\frac{\sqrt{\lambda_1 (\Omega^{-1}M^T W^{-1}M)}}{
1 + \sqrt{1+\lambda_1 (\Omega^{-1}M^T W^{-1}M)}}\right)^{2}
.\:  \label{2.7}
\eeq 
\end{lemma}
\vskip 2ex
Lemma \ref{lem1} is the key result we will use to prove that iteration 
(\ref{3.7}) converges for any $P_0\in\cal P$, and thus the risk sensitive like filter converges. 
 

\section{$N$-block risk sensitive like Riccati mapping}\label{sec:Nblock}
The robust filter (\ref{3.3})-(\ref{condizione_su_gamma}) can be interpreted as solving a standard
least-square  filtering problem with time-varying parameters in Krein space.
The Krein state space model consists of dynamics (\ref{3.1})
and observations (\ref{3.2}), to which we must adjoin the
new observations
\beq
0 = x_t + v_t^R \: . \label{4.6}
\eeq
The components of noise vectors $u_t$, $v_t$ and $v_t^R$
now belong to a Krein space and have the inner product
\beq
\left \langle \bmat {c}
u_t\\
v_t \\
v_t^R
\emat
\, , \,
\bmat {c}
u_s\\
v_s \\
v_s^R
\emat \right \rangle
= \bmat {ccc}
I_m & 0 &  0 \\
 0  & I_p &0 \\
 0 & 0 & -(\theta_{t-1})^{-1} I_n \\
\emat \delta_{t-s} .\label{4.7}
\eeq Note that, in the classical risk sensitive framework, \cite{HSK2,HSK3},
$v_t^{R}$ with $t\geq 0$ are identically distributed, whereas are not in this setting.
Since $x_t$ is
Gauss-Markov, the downsampled process $x_k^d = x_{kN}$, with
$N$ integer, is also Gauss-Markov with state space model
\bea
x_{k+1}^d &=& A^N x_k^d + {\cal R}_N {\bf u}_k^N \label{3.9} \\
{\bf y}_k^N &=& {\cal O}_N x_k^d +{\cal D}_N {\bf v}_k^N+ {\cal H}_N {\bf
u}_k^N \label{3.10} \\
{\bf 0} &=& {\cal O}_N^R x_k^d + {\bf v}_k^{R N} + {\cal L}_N {\bf u}_k^N
\:  \label{4.8}
 \eea
  where \bea
    {\bf u}_k^N &\triangleq & \bmat {cccc}
u_{kN+N-1}^T & u_{kN+N-2}^T & \ldots & u_{kN}^T \emat^T \nnl
{\bf v}_k^N & \triangleq & \bmat {cccc}
v_{kN+N-1}^T & v_{kN+N-2}^T & \ldots & v_{kN}^T \emat^T \nnl {\bf
y}_k^N & \triangleq& \bmat {cccc} y_{kN+N-1}^T & y_{kN+N-2}^T & \ldots &
y_{kN}^T \emat^T \nnl 
{\bf v}_k^{RN} & \triangleq& \left[\begin{array}{ccc}(v_{kN+N-1}^R)^T &
(v_{kN+N-2}^R)^T & \ldots \end{array}\right. \nnl
&& \hspace{3.5cm}\left. \begin{array}{cc}  \ldots & (v_{kN}^R)^T \end{array}\right] ^T .\nn \nn \eea In the model
(\ref{3.9})--(\ref{4.8}) \bea {\cal R}_N & \triangleq & \bmat {cccc} B & AB
& \ldots & A^{N-1}B \emat \nnl {\cal O}_N & \triangleq & \bmat {cccc}
(CA^{N-1})^T & \ldots & (CA)^T & C^T \emat^T \nn \\
{\cal O}_N^R & \triangleq&  \bmat {cccc}
(A^{N-1})^T & \ldots & (A)^T & I
\emat^T  \nn\\
{\cal D}_N & \triangleq & I_N \otimes D.
 \eea
Note that ${\cal R}_N$ and ${\cal O}_N$ denote
respectively the $N$-block reachability and observability matrices
of system (\ref{3.1})--(\ref{3.2}), where the blocks forming
${\cal O}_N$ are written from bottom to top instead of the usual
top to bottom convention. If the pairs $(A,B)$ and $(C,A)$ are
reachable and observable, ${\cal R}_N$ and ${\cal O}_N^T$ have full
row rank for $N \geq n$. In (\ref{3.10}) and (\ref{4.8}), if
\bea
H_t & \triangleq & \left \{ \begin{array} {cc}
CA^{t-1}B & t \geq 1 \\
0 & \mbox{otherwise }
\end{array} \right. \nn\\
  L_t & \triangleq & \left \{ \begin{array} {cc}
A^{t-1}B & t \geq 1 \\
0 & \mbox{otherwise}
\end{array} \right.\nn
\eea
${\cal H}_N$ and ${\cal L}_N$ are block
Hankel matrices defined as follows
\bea
{\cal H}_N & \triangleq & \bmat {cccccc}
0 & H_1 & H_2 & \cdots & H_{N-2} &H_{N-1}\\
0 & 0   & H_1 &  H_2 & \cdots & H_{N-2}\\
0 & 0   &  0 & H_1 &  \cdots & H_{N-3}\\
\vdots & \vdots &  \vdots &  & \vdots & \vdots\\
0 & 0   & \cdots & \cdots & 0 & H_1 \\
0 & 0   & \cdots & \cdots & \cdots & 0
\emat\nn\\
{\cal L}_N & \triangleq & \bmat {cccccc}
0 & L_1 & L_2 & \cdots & L_{N-2} &L_{N-1}\\
0 & 0   & L_1 &  L_2 & \cdots & L_{N-2}\\
0 & 0   &  0 &  L_1 &  \cdots & L_{N-3}\\
\vdots & \vdots &  \vdots &  & \vdots & \vdots\\
0 & 0   & 0 & \cdots & \cdots & L_1 \\
0 & 0   & 0 & \cdots & \cdots & 0\\
\emat  \nn.\eea

The Krein space inner product of the observation noise vector
\[
\bmat{c}
{\bf w}_k^N\\
{\bf w}_k^{R N}
\emat = \bmat {c}
{\cal D}_N{\bf v}_k^N\\
{\bf v}_k^{R N}
\emat + \bmat {c}
{\cal H}_N \\
{\cal L}_N
\emat {\bf u}_k^N
\]
admits the following decomposition
\bea
&&  {\cal K}_{\Theta_{N,k} } \defn \left \langle \bmat {c}
{\bf w}_k^N \\
{\bf w}_k^{R N}
\emat \, , \, \bmat {c}
{\bf w}_k^N \\
{\bf w}_k^{R N}
\emat \right \rangle \nnl
&=& \bmat {cc}
{\cal D}_N{\cal D}_N^T & 0 \\
0 & -(\Theta_{N,k})^{-1}
\emat + \bmat {c}
{\cal H}_N \\
{\cal L}_N \emat \bmat {cc}
{\cal H}_N^T & {\cal L}_N^T
\emat \nnl
&=& \bmat {cc}
I_{Np} & 0 \\
{\cal L}_N{\cal H}_N^T ({\cal D}_N{\cal D}_N^T + {\cal H}_N{\cal H}_N^T)^{-1} & I_{Nn}
\emat\nn\\
&& \hspace*{0.3cm} \times
\bmat {cc}
{\cal D}_N{\cal D}_N^T + {\cal H}_N{\cal H}_N^T & 0 \\
0 & S_{\Theta_{N,k}}
\emat \nnl
&& \hspace*{0.3cm} \times \bmat {cc}
I_{Np} & ({\cal D}_N{\cal D}_N^T + {\cal H}_N{\cal H}_N^T)^{-1} {\cal H}_N{\cal L}_N^T \\
0 & I_{Nn}
\emat \: , \label{4.9}
\eea
where
\beq
S_{\Theta_{N,k}} \defn - (\Theta_{N,k})^{-1}  + {\cal L}_N(I_{Nm}+{\cal H}_N^T ({\cal D}_N {\cal D}_N^T)^{-1}{\cal
H}_N)^{-1}{\cal L}_N^T \label{4.10}
\eeq
denotes the Schur complement of the $(1,1)$ block inside ${\cal
K}_{\Theta_{N,k}}$ with
\beq \Theta_{N,k}  \defn  \mathrm{diag}(\theta_{kN+N-2} ,\theta_{kN+N-3} ,\ldots , \theta_{kN-1} )\otimes I_n.\eeq
 The projection of noise vector ${\bf u}_k^N$ on
the Krein subspace spanned by the observation noise vector
$\bmat {cc} ({\bf w}_k^N)^T & ({\bf w}_k^{R N})^T \emat^T$ is
then given by
\[
\hat{\bf u}_k^N = \bmat {cc}
G_k &  G_k^R
\emat \bmat {c}
{\bf w}_k^N \\
{\bf w}_k^{N R}
\emat \: ,
\]
where
\[
\bmat {cc}
G_k & G_k^{R }
\emat = \bmat {cc}
{\cal H}_N^T & {\cal L}_N^T
\emat ({\cal K}_{\Theta_{N,k}})^{-1} \: ,
\]
and the residual $\tilde{\bf u}_k^N = {\bf u}_k^N -\hat{\bf u}_k^N$
has for inner product
\bea
&& {\cal Q}_{\Theta_{N,k}} \defn \langle \tilde{\bf u}_k^N \, , \, \tilde{\bf u}_k^N \rangle
\nn\\
&& \hspace{0.3cm}=I_{Nm} - \bmat {cc}
G_{k} & G_{k}^{ R }
\emat {\cal K}_{\Theta_{N,k}}
\bmat {c}
G_k^T \\ (G_k^{ R })^T \emat  \nnl[2ex]
&& \hspace{0.3cm}= [ I_{Nm} + {\cal H}_N^T ({\cal D}_N {\cal D}_N^T)^{-1} {\cal H}_N -  {\cal L}_N^T \Theta_{N,k}{\cal
L}_N ]^{-1} \nn\: .
\eea
Then by multiplying the observation equation obtained by
combining equations (\ref{3.10}) and (\ref{4.8}) by
${\cal R}_N \bmat {cc}  G_k &
 G_k^{R} \emat$ and subtracting it from
(\ref{3.9}), we obtain the state space equation
\beq
x_{k+1}^d = \alpha_{N,k} x_k^d + {\cal R}_N \tilde{\bf u}_k^N
+{\cal R}_N G_k {\bf y}_k^N
\label{4.13}
\eeq
with
\[
\alpha_{N,k} \defn A^N - {\cal R}_N [ G_k {\cal O}_N
+ G_k^{R } {\cal O}_N^R ] \: ,
\]
where the driving noise $\tilde{\bf u}^N_k$is now orthogonal to the noises ${\bf w}_k^N$
and ${\bf w}_k^{R N}$ appearing in observation equations (\ref{3.10})
and (\ref{4.8}). Accordingly, the time-varying Riccati iteration associated to the
downsampled model takes the form
\beq
P_{k+1}^d = r^d_c (P_k^d)   \label{4.14}
\eeq
where
\beq r^d_c (P_k^d) \defn \alpha_{N,k}
[ (P_k^d)^{-1} + \Omega_{\Theta_{N,k}} ]^{-1}  \alpha_{N,k}^T
+ W_{\Theta_{N,k}}  \nn \eeq
\bea
\Omega_{\Theta_{N,k}} & \triangleq & \bmat {cc}
{\cal O}_N^T & ({\cal O}_N^R)^T
\emat {\cal K}_{\Theta_{N,k}}^{-1} \bmat {c}
{\cal O}_N \\
{\cal O}_N^R
\emat \nnl
&=& \Omega_N + {\cal J}_N^T S_{\Theta_{N,k}}^{-1} {\cal J}_N
\label{4.15}
\eea
with
\bea
{\cal J}_N &\defn & {\cal O}_N^R  -{\cal L}_N{\cal H}_N^T [{\cal D}_N{\cal D}_N^T+{\cal H}_N{\cal
H}_N^T]^{-1}{\cal O}_N\nn\\
\Omega_N & \defn & {\cal O}_N^T({\cal D}_N{\cal D}_N^T+{\cal H}_N{\cal H}_N^T)^{-1} {\cal O}_N\nn
\eea
and
\beq
W_{\Theta_{N,k}} \triangleq {\cal R}_N{\cal Q}_{\Theta_{N,k}}{\cal R}_N^T
\:. \label{4.16}
\eeq

\section{Convergence of the risk sensitive like filter}\label{sec:contraction_property}

In this Section we show that the filter (\ref{3.3})-(\ref{condizione_su_gamma}) asymptotically  converges, or equivalently iteration (\ref{3.7}) converges for any $P_0\in\cal P$, for tolerance values sufficiently small.

The idea is to find an upper bound for $c$, say $c_{\mathrm{MAX}}$, 
for which the Gramians $\Omega_{\Theta_{N,k}}$ and $W_{\Theta_{N,k}}$ are positive definite for $k\geq \tilde q$ and $\tilde q$ is a fixed integer number.  Then, by Lemma \ref{lem1}, $r_c^d(\cdot)$ is a strict contraction for $k\geq \tilde q$. Since $r_c^d(\cdot)$ is the $N$-block mapping of $r_c^R(\cdot)$, we conclude that iteration (\ref{3.7}) converges.

Note that, $\Omega_{\Theta_{N,k}}$ and $W_{\Theta_{N,k}}$ depend on the positive definite matrix $\Theta_{N,k}$.
It is not difficult to see that, a sufficient condition to guarantee ${\cal Q}_{\Phi}$ positive definite for $0 \preceq  \Phi \prec \tilde \phi_N I_{Nn}$, and thus also $W_{\Phi}$ for $N\geq n$, is that
\[ \tilde \phi_N=\frac{1}{\lambda_1({\cal L}_N(I_{Nm}+ {\cal H}_N^T ({\cal D}_N{\cal D}_N^T)^{-1}{\cal H}_N)^{-1}{\cal L}_N^T)}>0.\]
Such a condition also guarantees that $S_{\Phi}^{-1}$ is negative definite. 
\begin{lemma}
\label{lem2}
Let $\Phi_1$ and $\Phi_2$ be such that $0\preceq \Phi_2\preceq \Phi_1\prec \tilde \phi_N I_{Nn}$. Then,
\bea  \Omega_{\Phi_1 } & \preceq & \Omega_{\Phi_2}\nn \\
W_{\Phi_1} & \succeq  & W_{\Phi_2}.
\eea
\end{lemma}
\vskip 2ex
Note that
\bea \left. \Omega_\Phi\right|_{\Phi=0}&=&\Omega_N\nn\\
\left. W_\Phi \right|_{\Phi=0} &=& {\cal R}_N ( I_{Nm}+{\cal H}_N^T ({\cal D}_N{\cal D}_N^T)^{-1}{\cal H}_N )^{-1}{\cal R}_N^T\nn\eea
which are positive definite matrices for $N\geq n$
because the pairs $(C,A)$ and $(A,B)$ are observable
and reachable, respectively. Accordingly, in view of Lemma \ref{lem2}, there exists a constant $\phi_N>0$ with $N\geq n$ such that
 \beq \label{pos_cond_gramians}\Omega_{\Phi}, W_{\Phi}\succ 0 ,\;\;\forall \,0 \preceq  \Phi \prec \phi_N I_{Nn}. \eeq
 As noticed before, $W_{\Phi}$ is positive definite over the range $0\preceq \Phi \prec \tilde \phi_n I_{Nn}$. Then,
we set $\phi_N=\tilde \phi_N$ and
check whether
$\Omega_{\phi_N I_{Nn}}$ is positive definite or not. If not, we decrease $\phi_N$ up to $\Omega_{\phi_N I_{Nn}}$ becomes positive semi-definite but singular. In this way, (\ref{pos_cond_gramians}) holds. Therefore, 
if $\Theta_{N,k}\prec \phi_N I_{Nn}$ for $k\geq \tilde q$, then the Gramians $\Omega_{\Theta_{N,k}}$ and $W_{\Theta_{N,k}}$ are  positive definite for $k\geq \tilde q$, and thus iteration (\ref{3.7}) converges.

Now, we characterize an upper bound for $c$, say $c_{MAX}$, which guarantees that $\Theta_{N,k} \prec \phi_N I_{Nn}$ for $k\geq \tilde q$. To this aim,
we need the following two Propositions.

\vskip 2ex \begin{proposition} \label{proposition_lower_bound}Consider iterations (\ref{3.7}) and (\ref{4.14}) with an arbitrary $P_0\in{\cal P}$. Moreover, consider the iteration 
 \beq \overline{P}_{t+1}=r(\overline{P}_t),\;\; \overline{P}_0=BB^T.\nn \eeq Then, after a finite number of steps, say $q+1$, we have \[ P_t\succeq \overline{P}_q,\;\; t\geq q+1, \] and after $\tilde q=\lceil \frac{q+1}{N}\rceil$ steps \[ P_k^d\succeq \overline{P}_q,\;\; k\geq \tilde q.\]\end{proposition} \vskip 2ex

\begin{proposition} \label{prop_gamma}
Assuming that $0<\theta< (\lambda_1(P))^{-1}$, the following facts hold:\begin{enumerate}
                           \item $\gamma(\cdot,P)$ is monotone increasing over $\mathbb{R}_+$
                           \item $\gamma(\theta,P)>0$ for any $P\in\bar{\cal P}$ with $P\neq 0$
                           \item If $ P\succeq Q$ then $\gamma(\theta,P)\geq \gamma(\theta,Q)$
                         \end{enumerate}
\end{proposition}
\vskip2ex
By Proposition \ref{proposition_lower_bound}, $P_t\succeq \overline{P}_q$ with $t\geq q+1$ and $q$ is fixed. Then, by Proposition
\ref{prop_gamma}, (\ref{condizione_su_gamma}) implies that \[ \theta_{t-1}\leq \bar \theta,\;\ \forall t\geq q+1\]
where $\bar \theta$ is such that $\gamma(\theta_{t-1},P_t)=\gamma(\bar \theta,\overline{P}_q)$. Thus, condition 
$\Theta_{N,k}< \phi_N I_{Nn}$ for $k\geq \tilde q$
 is guaranteed if we choose $c$ in such a way that $\bar \theta<\phi_N$. The idea is formalized in the following Proposition.  

\vskip2ex
\begin{proposition} \label{corollary} Let $c$ be such that $0< c<c_{MAX}$ with $c_{MAX}\triangleq\gamma(\phi_N,\overline{P}_q)$, $N\geq n$, and $q$ fixed. Then, the mapping $r^d_c(\cdot)$ is strictly contractive after $\lceil\frac{q+1}{N}\rceil$ steps. Accordingly, iteration 
(\ref{3.7}) converges to a unique solution for any initial condition $P_0\in{\cal P}$.
\end{proposition}
\vskip 2ex
Note that, by Proposition \ref{prop_gamma}, the map \[ q\mapsto \gamma(\phi_N,\overline{P}_q) \]
is nondecreasing. Thus, we have to choose $q$ sufficiently large in order to find a bigger $c_{MAX}$.

\section{An Example}\label{sec:example}
\begin{figure}[thpb]
      \centering
      \includegraphics[width=9.2cm]{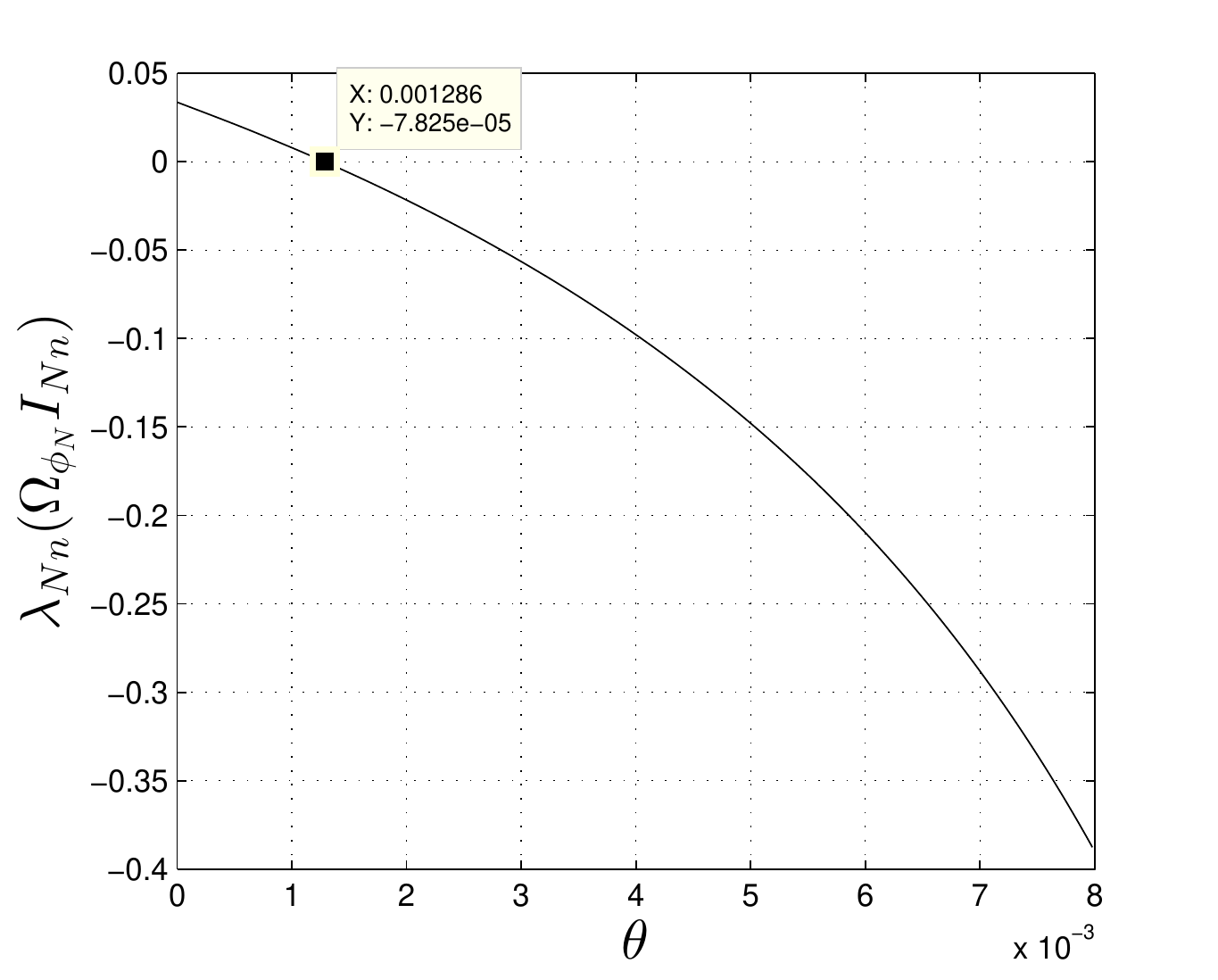}
      \caption{Minimum eigenvalue of $\Omega_{\phi_N I_{Nn}}$ over $\phi_N\in[0,8\cdot 10^{-3}]$.}
      \label{omin}
   \end{figure}

We consider the Gauss-Markov state space model earlier employed in \cite{LZ}
\bea \label{example_model} & A=\left[
       \begin{array}{cc}
         0.1 & 1 \\
         0 & 1.2 \\
       \end{array}
     \right], & B=I_2\nn\\
  & C=\left[
         \begin{array}{cc}
           1 & -1 \\
         \end{array}
       \right],
       & D=1
\eea with $n=2$, $m=2$ and $p=1$.
We select $N=8$, in this way $N\geq n$. Note that, larger values of $N$ can be considered. We find that
$\tilde \phi_N=1.6\cdot 10^{-2}$. In Figure \ref{omin}
 we depict the smallest eigenvalue of $\Omega_{\phi_N I_{Nn}}$ over the range
$\phi_N\in [0,8\cdot 10^{-3}]$. We find it becomes zero for $\phi_N \cong 1.3\times 10^{-3}$.

\begin{figure}[thpb]
      \centering
      \includegraphics[width=9.2cm]{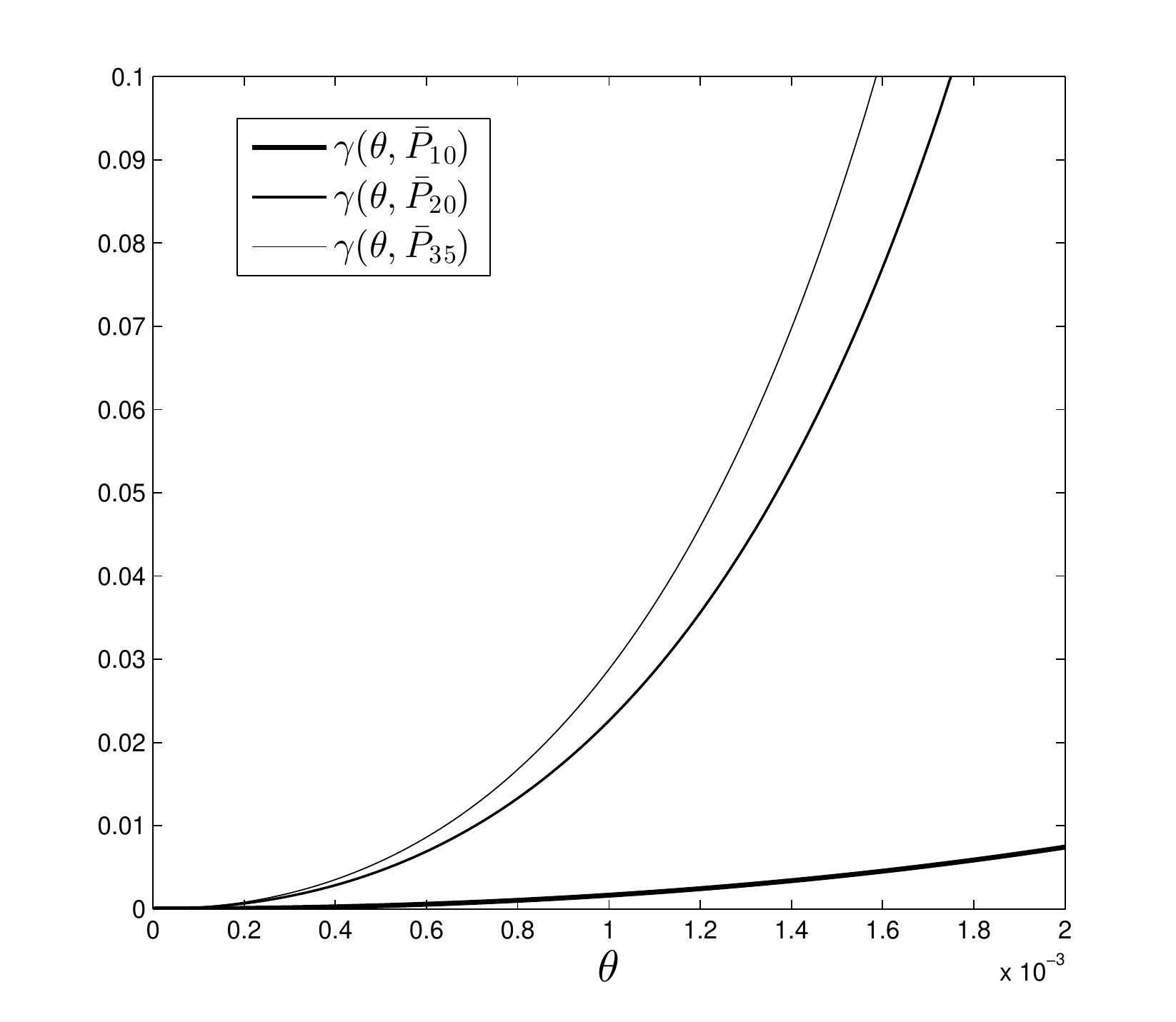}
      \caption{$\gamma(\theta,\bar{P}_q)$ with $q=10$, $q=20$ and $q=35$.}
      \label{gamma_plot}
   \end{figure}  In Figure \ref{gamma_plot} we depict $\gamma(\theta,\bar{P}_{10})$, $\gamma(\theta,\bar{P}_{20})$ and $\gamma(\theta,\bar{P}_{35})$. Note that,
 $\gamma(\phi_N,\overline{P}_{10})\cong2.9\cdot 10^{-3}$, $\gamma(\phi_N,\overline{P}_{20})\cong4.39\cdot 10^{-2}$  and $\gamma(\phi_N,\overline{P}_{35})\cong5.43\cdot 10^{-2}$. As expected, it is better to choose $\bar{P}_{35}$ for which we have $c_{\mathrm{MAX}}\cong5.43\cdot 10^{-2}$. We conclude that the risk sensitive like filter (\ref{3.3})-(\ref{condizione_su_gamma}) having tolerance parameter $c< 5.43\cdot 10^{-2}$
and nominal model (\ref{example_model}) asymptotically converges to a unique solution.

\section{Conclusion}
We analyzed the convergence of a risk sensitive like filter subject to an incremental tolerance. By contraction analysis, we showed that the corresponding $N$-block risk sensitive like Riccati mapping is strictly contractive for tolerance   values sufficiently small. Accordingly, the corresponding iteration converges to a fixed point, and thus the robust filter converges.

\appendix

\subsection{Proof of Lemma \ref{lem2}}

Let $\Phi$ and $\Delta$ belong to $\cal P$ and $\overline{{\cal P}}$, respectively, and such that $0\preceq \Phi+\Delta\prec \tilde \phi_n I_{Nn}$. The first variation of
$S_{\Phi}^{-1}$ and ${\cal Q}_{\Phi}$ with respect to $\Phi$ in direction $\Delta$
are, respectively,
 \[
\delta S_{\Phi; \Delta}^{-1}= -(S_{\Phi})^{-1}
\Phi^{-1}\Delta \Phi^{-1}
S_{\Phi}^{-1} \preceq 0
\]
and
\[
\delta {\cal Q}_{\Phi;\Delta}= {\cal Q}_{\Phi}
{\cal L}_N^T \Delta {\cal L}_N
 {\cal Q}_{\Phi}\succeq 0 \: .
\] This implies that \bea  S^{-1}_{\Phi_1} &\preceq & S^{-1}_{\Phi_2}\nn\\
   {\cal Q}_{\Phi_1} &\succeq & {\cal Q}_{\Phi_2}\nn \eea
   with $0\preceq \Phi_2\preceq \Phi_1\preceq \tilde \phi_N I_{Nn}$. Accordingly,
   \bea  \Omega_{\Phi_1} & =& \Omega_N +{\cal J}_N ^T S^{-1}_{\Phi_1}{\cal J}_N \preceq  \Omega_N+{\cal J}_N^T S^{-1}_{\Phi_2}{\cal J}_N =\Omega_{\Phi_2}\nn\\
  W_{\Phi_1}&=& {\cal R}_N {\cal Q}_{\Phi_1} {\cal R}_N^T  \succeq  {\cal R}_N {\cal Q}_{\Phi_2} {\cal R}_N ^T=W_{\Phi_2}.\nn \eea \qed\\

\subsection{Proof of Proposition \ref{proposition_lower_bound}}

Before proving the statement, we need the following two Lemmas. The first one regards
a risk sensitive mapping property, \cite[page 379]{HSK1}.
\vskip 2ex 
\begin{lemma}  Let $P\in{\cal P}$ such that $P^{-1}-\theta_1 I\in{\cal P}$ and $\theta_1\geq \theta_2\geq 0$. Then,
\[ r_{\theta_1}^{RS}(P)\succeq r_{\theta_2}^{RS}(P).\]
\end{lemma}
 \vskip 2ex
\begin{lemma} \label{prop_q}
Consider the sequence
\[ \tilde P_{t+1}=r^R_c(\tilde P_t),\;\; \tilde P_0\succeq BB^T. \]
Then, \beq \tilde P_t\succeq \overline{P}_q,\;\;\forall \; t\geq q. \eeq
\end{lemma}
\vskip 2ex
\begin{proof}
We prove by induction that
\beq \label{statement_lower_bound} \tilde P_t\succeq \overline{P}_t,\;\;\forall \; t\geq 0 .\eeq
Since, the sequence $\{ \overline{P}_t\}$ is nondecreasing, see for instance \cite{KSH}, then the statement follows for a fixed value of $q$.
 For $t=0$, we have $\tilde P_0\succeq \overline{P}_0$. Assume that (\ref{statement_lower_bound}) holds at time $t$, then
\beq \tilde P_{t+1}=r_c^R(\tilde P_t)=r_{\theta_{t-1}}^{RS}(\tilde P_t)\succeq r(\tilde P_t)\succeq  r(\overline{P}_t)=\overline{P}_{t+1}\eeq accordingly (\ref{statement_lower_bound}) also holds at time $t+1$.
\end{proof}
\vskip 2ex

Now, we proceed with the proof of Proposition \ref{proposition_lower_bound}.
 Consider the sequence generated by (\ref{3.7}) with an arbitrary initial condition $P_0\in{\cal P}$. Note that $r_c^R(P)\succeq BB^T$ for any $P\in{\cal P}$.
Accordingly $P_t\succeq BB^T$ for $t\geq 1$. We define the sequence $\{\tilde P_t\}$ with $\tilde P_t=P_{t+1}$ and $\tilde P_0=P_1\succeq BB^T$. In view of Lemma \ref{prop_q}, we have that $\tilde P_t\succeq \overline{P}_q$ after $q$ steps. This implies $ P_t\succeq \overline{P}_q$ for $t\geq q+1$. Noting that $P_k^d=P_{kN}$, the last statement follows.
\qed

\subsection{Proof of Proposition \ref{prop_gamma}}

 The first point has been proved in \cite{LN}. Regarding the second point, $\gamma(\theta,P)$ is equal to the information divergence among
the positive definite matrices $(I-\theta P)$ and $I$. Since $I-\theta P\neq I$, we get $\gamma(\theta,P)>0$.
 In order to prove the last point, we compute the first variation
of $\gamma(\theta,P)$ with respect to $P$ in direction $Q\in \bar{\cal P} $:
\bea && \delta \gamma(\theta,P; Q)\nn\\ &&=\hspace{0.2cm} \frac{\theta}{2}\mathrm{tr}[-(I-\theta P)^{-1}Q+(I-\theta P)^{-1}Q(I-\theta P)^{-1}]\nn\\
&&\hspace{0.2cm} =\frac{\theta}{2}\mathrm{tr}[Q^{\frac{1}{2}}(I-\theta P)^{-\frac{1}{2}}(-I+(I-\theta P)^{-1})\nn\\
&& \hspace{0.5cm}\times(I-\theta P)^{-\frac{1}{2}}Q^{\frac{1}{2}}]\succeq 0\nn\eea
where we exploited the fact that $(I-\theta P )^{-\frac{1}{2}}$ and $-I+(I-\theta P )^{-1}$ commute.\qed

\subsection{Proof of Proposition \ref{corollary}}

 Consider iterations (\ref{3.7}) and (\ref{4.14}). As showed in Proposition \ref{proposition_lower_bound}, after a finite number of steps, that is $q+1$ and $\tilde q=\lceil \frac{q+1}{N}\rceil$, we have, respectively, \[ P_t\succeq \overline{P}_q,\;\; t\geq q+1\]
\[ P^d_k\succeq \overline{P}_q,\;\; k\geq \tilde q.\]
Since $c<\gamma(\phi_N,\overline{P}_q)$, by Proposition \ref{prop_gamma} we have $\theta_{t}< \phi_N$ for $t\geq q+1$ and therefore $\Theta_{N,k}\prec \phi_N I_{Nn}$ for $k\geq \tilde q$.
Accordingly, the Gramians, $\Omega_{\Theta_{N,k}}$ and $W_{\Theta_{N,k}}$ are positive definite for $k\geq \tilde q$. By Lemma \ref{lem1}, the mapping $r^d_c(\cdot)$ is strictly contractive after $\tilde q$ steps.  Since $r^d_c(\cdot)$ is the $N$-block mapping of  $r_c^R(\cdot)$, it follows that the sequence generated by (\ref{3.7}) converges for any $P_0\in{\cal P}$.
\qed

\end{document}